\documentclass[preprint]{imsart}

\RequirePackage[OT1]{fontenc}
\RequirePackage{amsthm,amsmath}
\RequirePackage[]{natbib}
\RequirePackage[colorlinks,citecolor=blue,urlcolor=blue]{hyperref}
\usepackage{url}

\startlocaldefs
\numberwithin{equation}{section}
\theoremstyle{plain}
\newtheorem{theorem}{Theorem}[section]
\newtheorem{lemma}{Lemma}[section]
\newtheorem{proposition}{Proposition}[section]
\newtheorem{remark}{Remark}[section]
\newtheorem{example}{Example}[section]
\endlocaldefs

\begin{document}

\begin{frontmatter}
\title{Minimum Contamination and $\beta$-Aberration Criteria for Screening Quantitative Factors}
\runtitle{Contamination and $\beta$-Aberration}

\begin{aug}
\author{\fnms{Po} \snm{Yang}
}
\and
\author{\fnms{Chang-Yun} \snm{Lin}\thanksref{t1}
\ead[label=e2]{chlin6@nchu.edu.tw}
\ead[label=u1, url]{http://www.amath.nchu.edu.tw/\textasciitilde cylin/}
}

\thankstext{t1}{Correspondence should be addressed to Chang-Yun Lin. (e-mail: chlin6@nchu.edu.tw)}
\runauthor{Yang, P. and Lin, C.-Y.}

\affiliation{University of Manitoba, Canada\\ and\\ National Chung Hsing University, Taiwan}

\address{Department of Statistics, University of Manitoba, Canada\\
}

\address{Department of Applied Mathematics and Institute of Statistics\\ 
National Chung Hsing University, Taiwan\\
\printead{e2}\\
\printead{u1}
}
\end{aug}

\begin{abstract}

Tang and Xu [\textit{Biometrika} \textbf{101} (2014) 333--350] applied the minimum $\beta$-aberration criterion to selecting optimal designs for screening quantitative factors. They provided a statistical justification showing that minimum $\beta$-aberration criterion minimizes contamination of nonnegligible $k$th-order effects on the estimation of linear effects for $k=2,\cdots,r$, where $r$ is the strength of a design. Unfortunately, this result does not hold for $k>r$. In this paper, we provide a complete mathematical connection between $\beta$-wordlength patterns and contaminations (on the estimation of linear effects) and reveal that 
the minimum $\beta$-aberration criterion is not necessarily equivalent to the minimum contamination criterion for ranking designs.  
We prove that they are equivalent only when the number of factors of a design equals the strength plus one. We emphasize that the minimum $\beta$-aberration criterion, in fact, sequentially minimizes the contamination of nonnegligible $k$th-order effects on the estimation of the general mean, not on the estimation of linear effects. Therefore, the minimum contamination criterion should be more appropriate than the minimum $\beta$-aberration criterion for selecting optimal designs for screening quantitative factors.    
   
\end{abstract}

\begin{keyword}[class=MSC]
\kwd[Primary ]{62K05}
\kwd{62K15}
\end{keyword}

\begin{keyword}
\kwd{alias matrix}
\kwd{generalized minimum aberration}
\kwd{indicator function}
\kwd{mirror-symmetric design}
\end{keyword}

\end{frontmatter}

\section{Introduction}\label{sec:intro}
In design of experiments, it is important to comprehend the alias relationship between effects of fractional factorial designs. 
%
More serious aliases would lead to larger bias of estimations and affect statistical analysis results. 
For regular fractional factorial designs constructed on the Galois field, their alias structures can be easily observed through complete defining relations and are usually presented as {\it wordlength patterns}. 
%
Based on the hierarchy principle, the {\it minimum aberration criterion}, which ranks designs by sequentially minimizing wordlength patterns, is commonly used for selecting optimal designs. 
However, nonregular designs, such as the Plackett-Burman designs, do not have defining relations and have more complicated alias structures. The defining relation method does not work for this type of designs. 
A general method often introduced in text books or literatures is using the polynomial or regression model to generate the {\it alias matrix}, which captures the alias of specified model terms with terms that are potentially important but are not included in the model. 
The linear effect model is usually considered and the {\it contamination} (or {\it the overall alias}) of nonnegligible $k$th-order effects on the estimation of linear effects are measured by the square norm of the alias matrix. 
The {\it minimum contamination criterion} can be defined as sequentially minimizing contaminations for selecting optimal designs.
 



\citet{r11} 
developed the \textit{$J$-characteristics} to extend the concept of the wordlength pattern to two-level nonregular designs. Based on the $J$-characteristics, they proposed the \textit{minimum $G_2$-aberration criterion} for selecting optimal two-level regular or nonregular designs. 
\citet{r14} 
further extended the minimum $G_2$-aberration criterion and proposed the \textit{generalized minimum aberration criterion} for asymmetric fractional factorial designs. 
\citet{r11} 
and \citet{r14} 
investigated the connection between contaminations and aberrations and found that the minimum $G_2$-aberration criterion and the generalize minimum aberration criterion are equivalent to the minimum contamination criterion for ranking designs.    
These results are important because they provide statistical justification for the new proposed aberration criteria.  
%
%
%
%
%

Nevertheless, \citet{r1} 
pointed out that level permutation of a design could result in different geometric structures when factors are quantitative. 
They showed that the generalized minimum aberration criterion proposed in \citet{r14} 
can only distinguish combinatorially nonisomorphic designs for qualitative factors but not geometrically nonisomorphic designs for quantitative factors. 
To overcome this problem, \citet{r1} 
generalized the \textit{indicator function} proposed in \citet{r4} 
and \citet{r15} 
for designs with more than two levels by using the orthogonal polynomial basis. Based on the indicator function, they developed the \textit{$\beta$-wordlength pattern} to detect the geometrically structural change caused by level permutations and proposed the \textit{minimum $\beta$-aberration criterion} for selecting optimal designs when factors are quantitative. 

Although the minimum $\beta$-aberration criterion seems an intuitive extension of the generalized minimum aberration criterion and has been widely used for decade, its statistical justification was rarely discussed in literature. 
\citet{r12} 
studied the connection between $\beta$-wordlength patterns and contaminations and obtained an important theorem showing that minimum $\beta$-aberration criterion minimizes contamination of nonnegligible $k$th-order effects on the estimation of linear effects for $k=1,\cdots,r$, where $r$ is the strength of a design.
However, this theorem only provided a partial relationship between $\beta$-wordlength patterns and contaminations for $k\leq r$. The relationship of them for $k>r$ remains unsolved. Therefore, this result is not sufficient to provide statistical justification for the minimum $\beta$-aberration criterion. 
Meanwhile, we found that the minimum $\beta$-aberration design could be different from the minimum contamination design. A counterexample is given by two designs listed in Table~\ref{tb:2D}. 
We will later show in Example~\ref{ex:ineq} that design $D_1$ has less $\beta$-aberration while design $D_2$ has less contamination. 
This interesting result motivates us to investigate the complete relationship between the two criteria and to theoretically prove that they are inequivalent for ranking designs. 
We are also interested in examining their inconsistent ranking rates and exploring whether the two criteria could be equivalent under some conditions.  
We will discuss their statistical properties and rationalities for screening designs with quantitative factors.

%

The rest of this paper is organized as follows. Section~\ref{se:bk} introduces the indicator function and definitions of the minimum $\beta$-aberration and contamination criteria. Section~\ref{se:th} provides a mathematical equation showing the complete relationship between the two criteria. A condition in which the two criteria are equivalent is given. A comparison for two criteria and their statistical properties are discussed in Section~\ref{se:comp}. Section~\ref{se:re} is the concluding and remarks.

%







\section{Background and notation}\label{se:bk}

Let $D$ denote an orthogonal array with $n$ runs and $m$ factors $X_1,\cdots,X_m$, where the levels of factor $X_j$ are $0,1,\cdots,s_j-1$. Let $S_j=\{0,1,\cdots,s_j-1\}$. For factor $X_j$, let $c_0^j(x)=1$ and $c_u^j(x)$ be a polynomial of degree $u$ defined on $S_j$ for $u=1,\cdots,s_j-1$, such that 
\[
\sum_{x=0}^{s_j-1}c_u^j(x)c_v^j(x)=\left\{ 
\begin{array}{rl}
0,& u\neq v,\\
s_j & u=v.\\
\end{array}
\right.
\]
Let ${\cal T}=S_1\times \cdots\times S_m$ and $N=s_1\times\cdots\times s_m$. 
For a design point ${\bf x}=(x_1,\cdots,x_m)$, define 
\[
C_{\bf t}({\bf x})=\prod_{j=1}^mc_{t_j}^j(x_j),
\]
where ${\bf t}=t_1\cdots t_m\in {\cal T}$. The indicator function of $D$ is defined by a linear combination of $C_{\bf t}({\bf x})$ as
\[
F_{D}({\bf x})=\sum_{{\bf t}\in {\cal T}}b_{\bf t}C_{\bf t}({\bf x}),
\]
where the coefficient of $C_{\bf t}({\bf x})$ is uniquely determined by
\[
b_{\bf t}=\frac{1}{N}\sum_{{\bf x}\in D}C_{\bf t}({\bf x}). 
\] 
Define two norms: $||{\bf t}||_0=\sum_{j=1}^mI(t_j^+)$, where  $I(t_j^+)=1$ if $t_j=1,\cdots,s_j-1$ and 0 if $t_j=0$, which counts the number of nonzero elements in ${\bf t}$, and $||{\bf t}||_1=\sum_{j=1}^mt_j$, which calculates the polynomial degree of ${\bf t}$. 
For a design with strength $r$, coefficient $b_{\bf t}$ of the indicator function equals 0 if $0<||{\bf t}||_0\leq r$.
For quantitative factors, define the $\beta$-wordlength pattern by $(\beta_1,\cdots,\beta_{m'})$, where 
\[
\beta_k=\sum_{||{\bf t}||_1=k}\left(\frac{b_{\bf t}}{b_{\bf 0}}\right)^2
\]
for $k=1,\cdots,m'$ and $m'=\sum_{j=1}^m(s_j-1)$. The minimum $\beta$-aberration criterion is to sequentially minimize $\beta_k$ for $k=1,\cdots,m'$. The reader is referred to \citet{r1} 
for details.

Let ${\bf X}_{\bf t}=(C_{\bf t}({\bf x}_1),\cdots,C_{\bf t}({\bf x}_n))^T$ be an $n\times 1$ vector, where ${\bf x}_i=(x_{i1},\cdots,x_{im})$ is the $i$th run in $D$. 
For ${\bf t}\in{\cal T}$ with $||{\bf t}||_1=k$, ${\bf X}_{\bf t}$ is called the orthonormal polynomial contrast coefficient for the $k$th-order effect of the interaction $X_1^{t_1}\cdots X_m^{t_m}$. Then the full regression model for data analysis can be expressed by 
\[
{\bf Y}={\bf Z}_0\gamma_0+\sum_{k=1}^{m'}{\bf Z}_k\boldsymbol\gamma_k+\boldsymbol\epsilon,
\]
where ${\bf Y}$ is the $n\times 1$ vector of responses, $\gamma_0$ is the general mean and ${\bf Z}_0$ is an $n\times 1$ vector of 1s, $\boldsymbol\gamma_k$ is the vector of all $k$th-order effects and ${\bf Z}_k$ is the matrix of orthonormal polynomial contrast coefficients for $\boldsymbol\gamma_k$, and $\boldsymbol\epsilon$ is the $n\times 1$ vector of independent random errors. 
Now consider the common situation for screening experiments in which the linear effects are of primary interest and the fitted model is
\begin{equation}\label{eq:1st}
{\bf Y}={\bf Z}_0\gamma_0+{\bf Z}_1\boldsymbol \gamma_1+\boldsymbol\epsilon.
\end{equation}
%
%
The estimate of $\boldsymbol\gamma_1$ in model (\ref{eq:1st}) is $\hat  {\boldsymbol \gamma}_1=({\bf Z}_1^T{\bf Z}_1)^{-1}{\bf Z}_1^T{\bf Y}$.
Under the full model,
\[
E(\hat {\boldsymbol \gamma}_1)=\boldsymbol \gamma_1+\sum_{k=2}^{m'}{\bf A}_k\boldsymbol \gamma_{k},
\]
where
${\bf A}_k=({\bf Z}_1^T{\bf Z}_1)^{-1}{\bf Z}_1^T{\bf Z}_k$ is called the alias matrix.  
Define the contamination of non-negligible $k$th-order effects on the estimation of linear effects by 
\[
\lambda_k=||{\bf A}_k||^2=tr({\bf A}_k^T{\bf A}_k) 
\]
%
and the contamination pattern by $(\lambda_2,\cdots,\lambda_{m'})$. 
The minimum contamination criterion is to sequentially minimize $\lambda_k$ for $k=2,\cdots,{m'}$.
%

\begin{table}
\begin{center}
\caption{Two 18-run and 4-factor orthogonal arrays with 3 levels}
										
\begin{tabular}{c}
\hline		
\label{tb:2D}										
$\begin{array}{cc}																																					
&\\																																					
D_1&																																					
\begin{array}{c|cccccccccccccccccc}																																					
X_1	&	2	&	0	&	1	&	2	&	0	&	1	&	2	&	0	&	1	&	2	&	0	&	1	&	2	&	0	&	1	&	2	&	0	&	1	\\
X_2	&	0	&	1	&	2	&	0	&	1	&	2	&	1	&	2	&	0	&	2	&	0	&	1	&	1	&	2	&	0	&	2	&	0	&	1	\\
X_3	&	0	&	1	&	2	&	1	&	2	&	0	&	0	&	1	&	2	&	2	&	0	&	1	&	2	&	0	&	1	&	1	&	2	&	0	\\
X_4	&	0	&	1	&	2	&	1	&	2	&	0	&	2	&	0	&	1	&	1	&	2	&	0	&	0	&	1	&	2	&	2	&	0	&	1	\\
\end{array}\\																																					
&\\																																					
D_2&																																					
\begin{array}{c|cccccccccccccccccc}																																					
X_1	&	0	&	0	&	0	&	1	&	1	&	1	&	2	&	2	&	2	&	0	&	0	&	0	&	1	&	1	&	1	&	2	&	2	&	2	\\
X_2	&	2	&	0	&	1	&	2	&	0	&	1	&	2	&	0	&	1	&	2	&	0	&	1	&	2	&	0	&	1	&	2	&	0	&	1	\\
X_3	&	1	&	2	&	0	&	1	&	2	&	0	&	2	&	0	&	1	&	0	&	1	&	2	&	2	&	0	&	1	&	0	&	1	&	2	\\
X_4	&	0	&	1	&	2	&	2	&	0	&	1	&	1	&	2	&	0	&	1	&	2	&	0	&	2	&	0	&	1	&	0	&	1	&	2	\\
\end{array}\\																																					
&\\																																					
\end{array}$\\																					
\hline																
\end{tabular}																																																																									
\end{center}																			
\end{table}

\begin{example}\label{ex:ineq}
\normalfont
Designs $D_1$ and $D_2$ in Table~\ref{tb:2D} are two three-level orthogonal arrays with four factors and 18 runs. The polynomials for the $j$th factor with levels 0, 1, and 2 are $c_0^j(x)=1$, $c_1^j(x)={\sqrt \frac{3}{2}}(x-1)$, and $c_2^j(x)=\sqrt 2(\frac{3}{2}(x-1)^2-1)$. 
The $\beta$-wordlength patterns are (0, 0, 0.281, 0.797, 1.406, 0.313, 0.563, 0.141) for $D_1$ and (0, 0, 0.281, 0.844, 1.406, 0.781, 0.188, 0) for $D_2$.
The contamination patterns are (0.844, 2.203, 4.078, 2.109, 3.797, 0.688, 0.281) for $D_1$ and (0.844, 2.203, 3.984, 3.141, 2.953, 0.781, 0.094) for $D_2$.
According to the minimum $\beta$-aberration criterion, $D_1$ is a better design. However, $D_2$ is considered better than $D_1$ if the minimum contamination criterion is applied. 
\end{example}

\section{Connection between the two criteria for three-level designs}\label{se:th}

Example~\ref{ex:ineq} shows that the minimum $\beta$-aberration criterion and the minimum contamination criterion are not consistent for ranking designs. 
In this section, we investigate the mathematical connection between the two criteria through the indicator function for the commonly used $OA(n,3^m,r)$, the $n$-run and $m$-factor orthogonal array with 3 levels and strength $r$, where $r\geq 2$.  
%
Let ${\cal T}_{i,j}=\{{\bf t}\in{\cal T} \left|\right. ||{\bf t}||_0=i+j, ||{\bf t}||_1=i+2j\}$. Define 
\begin{equation}\label{eq:b.ij}
\beta_{i,j}=\sum_{{\bf t}\in{\cal T}_{i,j}}\left(\frac{b_{\bf t}}{b_{\bf 0}}\right)^2
\end{equation}
for $i=0,\cdots,m$, $j=0,\cdots,m-i$, and $\beta_{i,j}=0$, otherwise. 
Note that $i$ and $j$ represent the numbers of $1$ and $2$ in $\bf t$, respectively. 
%
%
Let ${\bf Z}_{i,j}$ be an $n\times {{m\choose {i}} {{m-i}\choose j}}$ matrix whose columns are the orthonormal polynomial contrast coefficients ${\bf X}_{\bf t}$s with $||{\bf t}||_0=i+j$ and $||{\bf t}||_1=i+2j$. Let ${\bf A}_{i,j}=({\bf Z}_1^T{\bf Z}_1)^{-1}{\bf Z}_1^T{\bf Z}_{i,j}=n^{-1}{\bf Z}_1^T{\bf Z}_{i,j}$ and define
\begin{equation}
\lambda_{i,j}=||{\bf A}_{i,j}||^2=tr({\bf A}_{i,j}^T{\bf A}_{i,j})
\end{equation}
for $i=0,\cdots,m$, $j=0,\cdots,m-i$, and $\lambda_{i,j}=0$, otherwise.
Let
\[
\xi_{i,j}=\sum_{l=1}^m\sum_{{\bf t}\in{\cal T}_{i,j}^{l(-)}}\prod_{g=0,2}\left(\frac{b_{{\bf t}|_{t_l=g}}}{b_{\bf 0}}\right)
\]
for $i=0,\cdots,m-1, j=0,\cdots,m-1-i$, and $\xi_{i,j}=0$, otherwise, 
where ${\cal T}_{i,j}^{l(-)}=\{{\bf t}\in{\cal T}\left| \right. ||{\bf t}||_0-I(t_l^+)=i+j,||{\bf t}||_1-t_l=i+2j\}$ and ${\bf t}|_{t_l=g}$ denote {\bf t} with $t_l=g$.   

\begin{lemma}\label{le:l.pq} 
For given $p$ and $q$, where $p+2q\geq 2$,
\begin{equation}\label{eq:l.pq}
\begin{array}{rl}
\lambda_{p,q}=&(p+1)\beta_{p+1,q}+\frac{p+1}{2}\beta_{p+1,q-1}+\frac{q+1}{2}\beta_{p-1,q+1}\\
&+(m-p-q+1)\beta_{p-1,q}+\sqrt 2\xi_{p-1,q}.
\end{array}
\end{equation}
\end{lemma}
%
%
%
%
Note that $\beta_{i,j}=0$ if $i,j<0$, $i,j>m$, or $i+j>m$, and $\xi_{i,j}=0$ if $i,j<0$, $i,j>m-1$, or $i+j>m-1$. Hence, some terms in equation (\ref{eq:l.pq}) may be zero. For instance, when $m=4$, $\lambda_{2,2}=\frac{3}{2}\beta_{3,1}+\frac{3}{2}\beta_{1,3}+\beta_{1,2}+{\sqrt 2}\xi_{1,2}$, $\lambda_{0,3}=\beta_{1,3}+\frac{1}{2}\beta_{1,2}$, and $\lambda_{4,0}=\frac{1}{2}\beta_{3,1}+\beta_{3,0}+\sqrt 2\xi_{3,0}$. For $p+q=m+1$, $\lambda_{p,q}=(m-p-q+1)\beta_{p-1,q}=0$ since $m-p-q+1=0$. Of course, $\lambda_{p,q}=0$ for $p+q>m+1$. 

Lemma \ref{le:l.pq} shows that $\lambda_{i,j}$ is a linear combination of $\beta_{i,j}$ and $\xi_{i,j}$. This mathematical result is powerful for developing the connection between minimum contamination and $\beta$-aberration criteria. 
Now let us express the $\beta$-wordlength pattern and the contamination pattern by
\begin{equation} \label{eq:b.k=b.pq}
\beta_k=\sum_{j=0}^{\lfloor k/2\rfloor}\beta_{k-2j,j},\mbox{ for }k=1,\cdots,m'
\end{equation}
and
\begin{equation}\label{eq:l.k=l.pq}
\lambda_k=\sum_{j=0}^{\lfloor k/2\rfloor}\lambda_{k-2j,j},\mbox{ for }k=2,\cdots,m',
\end{equation}
respectively, where $\lfloor w \rfloor$ is the largest integer not greater than $w$. 
Note that the summations of $\beta_{k-2j,j}$ and $\lambda_{k-2j,j}$ in equations (\ref{eq:b.k=b.pq}) and (\ref{eq:l.k=l.pq}) should start from $j=k-m$ when $k>m$ since $i+j$, where $i=k-2j$, must be smaller than or equal to the number of factors $m$. 
However, equations (\ref{eq:b.k=b.pq}) and (\ref{eq:l.k=l.pq}) still hold because $\beta_{i,j}=0$ and $\lambda_{i,j}=0$ for $i+j>m$ according to the definition.
Combining Lemma \ref{le:l.pq} and equations (\ref{eq:b.k=b.pq}) and (\ref{eq:l.k=l.pq}), we obtain the complete relationship between the contamination and the $\beta$-wordlength pattern as follows.

\begin{proposition} \label{pr:l.k}
For $k=2,\cdots,m'$, 
\begin{equation}\label{eq:l.k}
\lambda_k=(1+k-\frac{3}{2}\times\lceil \frac{k}{2}\rceil)\beta_{k+1}+(m-\frac{k-1}{2})\beta_{k-1}+B,
\end{equation}
where $B=\frac{3}{2}\sum_{j=0}^{\lceil k/2\rceil-1}(\lceil\frac{k}{2}\rceil-j)\beta_{k-2j+1,j}+\sqrt 2 \sum_{j=0}^{\lceil k/2\rceil-1}\xi_{k-2j-1,j}$; $\lceil w \rceil$ is the smallest integer not less than $w$.
\end{proposition}

It is not difficult to derive that equation (\ref{eq:l.k}) reduces to $\lambda_k=0$ for $k=2,\cdots,r-1$ and $\lambda_k=(k+1)\beta_{k+1}$ for $k=r$, which is consistent with Lemma 1 in \citet{r12} 
for three-level designs. 
However, Lemma~1 in \citet{r12} 
only gave a partial relationship between the contamination and the $\beta$-wordlength pattern for $k=2,\cdots,r$.  
The result presented in Proposition \ref{pr:l.k} provides a complete relationship for $k=2,\cdots,m'$. 
Equation~(\ref{eq:l.k}) shows that the contamination is not just a simple linear combination of the $\beta$-wordlength patterns. This result is important because it discloses that the minimum $\beta$-aberration criterion is not theoretically equivalent to the minimum contamination criterion for ranking designs. 
Therefore, the optimal designs obtained by the two criteria could be different. 

Although the relationship between the two criteria is not simple as shown in equation (\ref{eq:l.k}), there exists a condition in which the contamination can be expressed by a linear combination of $\beta$-wordlength patterns. 

\begin{proposition}\label{pr:r+1}
For $OA(n,3^m,r)$, when $m=r+1$,  
\begin{equation}\label{eq:r+1}
\lambda_k=\rho\beta_{k+1}+(r+1-\frac{k-1}{2})\beta_{k-1},
\end{equation}
where 
\begin{equation}\label{eq:rho}
\rho=\left\{
\begin{array}{ll}
k+1&\mbox{, if }k=2,\cdots, r,\\
\frac{1}{2}(3r+2-k)&\mbox{, if }k=r+1,\cdots, m'-1,\\
0&\mbox{, if }k=m'.
\end{array}
\right.
\end{equation}
\end{proposition}

%

Equation (\ref{eq:r+1}) is a general form for $k=2,\cdots,m'$. When $k\leq r+1$, it can be simplified as follows. For $k=2,\cdots,r-1$ (if $r\geq 3$), $\lambda_k=0$ since $\beta_{k+1}=0$ and $\beta_{k-1}=0$; for $k=r$, $\lambda_k=(k+1)\beta_{k+1}$ since $\beta_{k-1}=0$; for $k=r+1$, $\lambda_k=\frac{1}{2}(3r+2-k)\beta_{k+1}$ since $\beta_{k-1}=0$. 
%
Proposition~\ref{pr:r+1} shows that, when $m=r+1$, the minimum $\beta$-aberration criterion minimizes contamination of nonnegligible $k$th-order effects on the estimation of linear effects for $k=2\cdots,m'$.  
For instance, when $r=2$ and $m=3$, we obtain $\lambda_2=3\beta_3$, $\lambda_3=2\frac{1}{2}\beta_4$, $\lambda_4=2\beta_5+\frac{3}{2}\beta_3$, $\lambda_5=1\frac{1}{2}\beta_6+\beta_4$, and $\lambda_6=\frac{1}{2}\beta_5$. 
It is obvious that sequentially minimizing $\lambda_2$ to $\lambda_6$ is equivalent to sequentially minimizing $\beta_3$ to $\beta_6$. 
Therefore, the results of ranking designs by the two criteria are completely consistent when $m=r+1$.
We summarize this important result as the following theorem.

\begin{theorem}\label{th:eq}
For $OA(n,3^{r+1},r)$, the minimum $\beta$-aberration criterion is equivalent to the minimum contamination criterion.
\end{theorem}


There exists an interesting link between the two criteria for mirror-symmetric designs. For a given design, if we reverse the level order for all factors, we obtain its mirror-image. A design is called mirror-symmetric if its mirror-image is itself. \citet{r12} 
provided a theorem showing that a design is mirror-symmetric if and only if $\beta_k=0$ for odd $k$. Because $\lambda_k$ and $\beta_k$ have the relationship as shown in Proposition~\ref{pr:l.k}, we obtain a similar result for the minimum contamination criterion.  

\begin{theorem}\label{th:even}
A design is mirror-symmetric if and only if $\lambda_k=0$ for even $k$.
\end{theorem}

Theorem~\ref{th:even} implies that there is no alias between linear effects and nonnegligible $k$th-order effects for even $k$ when the design is mirror-symmetric.  
For a given design, we may eliminate the contamination of nonnegligible 2nd-order effects on the estimation of linear effects if the design is possible to become a mirror-symmetric design after level permutations. The following example demonstrates this interesting result.   

\begin{example}
\normalfont
For designs $D_1$ and $D_2$ in Table~\ref{tb:2D}, let us permute factor $X_2$ of $D_1$ and factor $X_3$ of $D_2$ by $\{0,1,2\}\rightarrow\{2,0,1\}$ to obtain $D_1'$ and $D_2'$, respectively. One can check that both $D_1'$ and $D_2'$ are mirror-symmetric designs. The contaminations for $D_1'$ and $D_2'$ are $(0,6,0,4.5,0,3.5,0)$ and $(0,5.063,0,7.313,0,1.625,0)$, respectively. This permutation successfully elimates the contamination of nonnegligible 2nd-order effects on the estimation of linear effects for both designs.




\end{example}





\section{Comparison and properties}\label{se:comp}

\begin{table}

\caption{Comparison of the minimum contamination criterion and the minimum $\beta$-aberration criterion for designs obtained from $L_{18}$.}

\begin{center}
\begin{tabular}{ccccc}									
\hline\label{tb:sum}									
$m$	&	Columns	&	$(\lambda_2,\lambda_3,\lambda_4)^{Rank}$	&	$(\beta_3,\beta_4,\beta_5)^{Rank}$	&	Consistent rate	\\
\hline									
3	&	$\breve{1},2,5$	&	$(0.000, 0.313, 1.500)^1$	&	$(0.000, 0.125, 0.750)^1$	&	$8/8$ $(100\%)$	\\
4	&	$\breve{1},2,3,\acute{5}$	&	$(0.000, 5.063, 0.000)^1$	&	$(1.875, 0.000, 1.625)^1$	&	$17/19$ $(89.5\%)$	\\
	&	$\breve{1},2,\acute{3},\acute{5}$	&	$(0.844, 2.203, 3.984)^9$	&	$(0.281, 0.844, 1.406)^{10}$	&		\\
	&	$2,3,4,\acute{5}$	&	$(0.844, 2.203, 4.078)^{10}$	&	$(0.281, 0.797, 1.406)^9$	&		\\
5	&	$\breve{1},2,3,\acute{4},\acute{5}$	&	$(0.000, 16.750, 0.000)^1$	&	$(0.000, 6.063, 0.000)^1$	&	$39/44$ $(88.6\%)$	\\
	&	$1,2,\acute{3},\breve{4},5$	&	$(2.531, 7.609, 12.070)^{21}$	&	$(0.844, 2.969, 3.492)^{22}$	&		\\
	&	$1,2,5,6,\acute{7}$	&	$(2.531, 7.734, 12.258)^{22}$	&	$(0.844, 3.094, 3.867)^{24}$	&		\\
	&	$\breve{1},2,3,4,\acute{5}$	&	$(2.531, 7.984, 11.766)^{24}$	&	$(0.844, 2.922, 3.750)^{21}$	&		\\
	&	$\breve{1},2,\acute{3},4,\acute{5}$	&	$(2.813, 5.969, 14.250)^{27}$	&	$(0.938, 2.313, 4.406)^{28}$	&		\\
	&	$1,2,3,5,6$	&	$(2.813, 6.094, 13.148)^{28}$	&	$(0.938, 2.297, 3.961)^{27}$	&		\\
6	&	$2,3,4,\breve{5},\breve{6},\breve{7}$	&	$(2.250, 19.875, 28.125)^1$	&	$(0.750, 6.938, 6.750)^1$	&	$38/40$ $(95.0\%)$	\\
	&	$1,2,3,5,6,\breve{7}$	&	$(5.063, 16.031, 29.391)^{24}$	&	$(1.688, 5.906, 7.828)^{25}$	&		\\
	&	$1,2,3,4,\breve{5},\breve{6}$	&	$(5.063, 16.031, 29.719)^{25}$	&	$(1.688, 5.766, 7.969)^{24}$	&		\\
7	&	$1,2,3,\breve{4},5,\breve{6},\breve{7}$	&	$(4.500, 41.063, 48.375)^1$	&	$(1.500, 14.625, 12.000)^1$	&	$29/29$ $(100\%)$	\\
\hline																	
\end{tabular}																									
		
\end{center}
\end{table}

We apply the two criteria to ranking $OA(18,3^m,2)$ for $m=3,\cdots,7$, where the $m$ columns are chosen from $L_{18}$ listed in Table 2 in \citet{r1}.  
Let $\acute{i}$ and $\breve{i}$ denote the level permutations for  column $i$ of $L_{18}$ by $\{0,1,2\}\rightarrow\{1,2,0\}$ and $\{0,1,2\}\rightarrow\{2,0,1\}$, respectively. The two permutations may obtain designs that have different geometric structures from the original design. 
We perform level permutations and calculate the contamination pattern and the $\beta$-wordlength pattern for each design. Only one design is kept for those having the same contamination pattern or $\beta$-wordlength pattern and ranks of the design are assigned according to the two criteria. The results are shown in Table~\ref{tb:sum}.
To save space, we only list optimal designs (with rank 1) and designs with inconsistent ranks assigned by the two criteria. 
In Table~\ref{tb:sum}, the first column is the number of factors, the second column lists the designs whose columns are chosen from $L_{18}$, the third column lists the first three contaminations and the superscript shows the rank assigned by the minimum contamination criterion, the fourth column lists the $\beta$-wordlength pattern for $k=3,4,5$ and the superscript shows the rank assigned by the minimum $\beta$-aberration criterion, and the last column calculates the consistent ranking rate between the two criteria, where the denominator is the total number of design ranks and the numerator is the total number of designs having consistent ranks assigned by the two criteria. 
For instance, when $m=4$, there are total 19 designs with different ranks by the two criteria. Among the 19 designs, designs with columns $1,2,\acute{3},\acute{5}$ and columns $2,3,4,\acute{5}$ are ranked as 9 and 10, respectively, according to the minimum contamination criterion, but are ranked as 10 and 9, respectively, by the minimum $\beta$-aberration criterion. The consistent ranking rate between the two criteria is $(19-2)/19=89.5\%$.         

Table~\ref{tb:sum} shows that the two criteria are not consistent for ranking designs when $m=4,5,6$. The lowest consistent rate ($88.6\%$) happens at $m=5$. When $m=3$, ranking designs by the two criteria are consistent, which verifies the result in Theorem~\ref{th:eq}. Although the two criteria are not theoretically equivalent when $m>r+1$, ranking designs by them may be completely consistent (e.g. $m=7$). It also shows that the optimal designs obtained by the two criteria are all identical in these cases. Observing the contamination patterns and the $\beta$-wordlength patterns, one can find that the optimal designs for $m=4,5$ are mirror-symmetric designs, which have $\beta_k=0$ for odd $k$ and $\lambda_k=0$ for even $k$ as the result in Theorem~\ref{th:even}. 
It is noteworthy to point out that, although we could permute levels for some designs to obtain mirror-symmetric designs, 
the mirror-symmetric designs 
may not be optimal designs. 
For instance, when $m=3$, there are three mirror-symmetric designs, columns $1,2,\acute{3}$, columns $\breve{1},2,\acute{5}$, and columns $1,3,\breve{4}$, ranked as 2, 3, and 4, respectively, by both criteria, but the optimal design shown in Table~\ref{tb:sum} 
is columns $\breve{1},2,5$, which is not a mirror-symmetric design.

In the view of data analysis, the minimum $\beta$-aberration criterion can be considered as a special case of the minimum contamination criterion. Let us replace the fitted mode (\ref{eq:1st}) by the general mean model
\begin{equation}\label{eq:mean}
{\bf Y}={\bf Z}_0\gamma_0+\boldsymbol\epsilon.
\end{equation}
The estimate of $\gamma_0$ is $\hat  {\gamma}_0=n^{-1}{\bf Z}_0^T{\bf Y}$.
Under the full model,
$
E(\hat {\gamma}_0)=\gamma_0+\sum_{k=1}^{m'}{\bf A}_k'\boldsymbol \gamma_{k},
$
where
$
{\bf A}_k'=n^{-1}{\bf Z}_0^T{\bf Z}_k
$ is a $1\times n_k$ vector whose elements are $(b_{\bf t}/b_{\bf 0})$s with $||{\bf t}||_1=k$, and $n_k$ is the number of $\bf t$s such that $||{\bf t}||_1=k$. The contamination of non-negligible $k$th-order effects on the estimation of the general mean is then defined by
\begin{equation}\label{eq:l.k'}
\lambda_k'=tr({\bf A}_k'^T{\bf A}_k')=\sum_{||{\bf t}||_1=k}(b_{\bf t}/b_{\bf 0})^2=\beta_k
\end{equation}
for $k=1,\cdots,m'$. That is, $\beta_k$ is identical to the contamination of non negligible $k$th-order effects on  the estimation of the general mean for $k=1,\cdots,m'$. 
In the screening experiments, we would rather fit a linear effect model than just fit the general mean model. Therefore, sequentially minimizing $\lambda_k$ should be more reasonable than sequentially minimizing $\beta_k$ ($\lambda_k'$) for screening quantitative factors. We summarize this result as the following theorem.

\begin{theorem}
The minimum $\beta$-aberration criterion minimizes contamination of nonnegligible $k$th-order effects on the estimation of the general mean for $k=1,\cdots,m'$.
\end{theorem}     


In practice, higher-order effects are difficult to explain and hence are usually considered negligible.  
For $OA(n,3^m,r)$, if $k$th-order effects are assumed negligible for $k>r$, then the full model reduces to   
\begin{equation}\label{eq:rd}
{\bf Y}={\bf Z}_0\gamma_0+\sum_{k=1}^{r}{\bf Z}_k\boldsymbol\gamma_k+\boldsymbol\epsilon.
\end{equation}
Under model (\ref{eq:rd}), the contamination pattern reduces to $(\lambda_2,\cdots,\lambda_{r-1},\lambda_r)=(0,\cdots,0,(r+1)\beta_{r+1})$. 
In fact, $\beta_k=0$ for $k=1,\cdots,r$. 
Therefore, minimum $\beta$-aberration criterion minimizes contamination of nonnegligible $k$th-order effects on the estimation of linear effects for $k=2,\cdots,r$. 
Based on this assumption, the minimum $\beta$-aberration criterion can be considered, technically, equivalent to the minimum contamination criterion. We summarize this result as the following theorem.

\begin{theorem}
If $k$th-order effects are assumed negligible for $k>r$, then the minimum contamination criterion is minimizing $\lambda_r$, which is equivalent to minimizing $\beta_{r+1}$.
\end{theorem}

Compare to the minimum $\beta$-aberration criterion, the minimum contamination criterion is more flexible for different situations. 
For example, when the variance-covariance matrix of error terms is $\boldsymbol\Sigma\neq\sigma^2\bf I$ (e.g. split-plot designs), the generalized least square estimate of $\boldsymbol\gamma_1$ in equation (\ref{eq:1st}) is 
\[
\hat {\boldsymbol\gamma}_1=({\bf Z}_1^T\boldsymbol\Sigma^{-1}{\bf Z}_1)^{-1}{\bf Z}_1^T\boldsymbol\Sigma^{-1}{\bf Y}
\] 
and 
$
E(\hat {\boldsymbol\gamma}_1)=\boldsymbol \gamma_1+\sum_{k=2}^{m'}{\bf A}_k^*\boldsymbol \gamma_{k},
$
where
\[
{\bf A}_k^*=({\bf Z}_1^T\boldsymbol\Sigma^{-1}{\bf Z}_1)^{-1}{\bf Z}_1^T\boldsymbol\Sigma^{-1}{\bf Z}_k
\]
for $k=2,\cdots,m'$. The minimum contamination criterion can be modified as sequentially minimizing 
\begin{equation}\label{eq:l.k*}
\lambda_k^*=tr({\bf A}_k^{*T}{\bf A}_k^*)
\end{equation} 
for $k=2,\cdots,m'$. We summarize the discussion as the following remark.

\begin{remark}
For multistratum designs with variate-covariate matrix ${\boldsymbol\Sigma}\neq\sigma^2{\bf I}$, the contamination of $k$th-order effects on the estimation of linear effects is $\lambda_k^*$ in (\ref{eq:l.k*}). The minimum contamination criterion is to sequentially minimize $\lambda_k^*$ for $k=2,\cdots,m'$.
\end{remark}






\section{Concluding and remarks}\label{se:re}

The minimum $\beta$-aberration criterion has been proposed and widely used for decade but its statistical justification was rarely discussed in literature. 
Although \citet{r12} attempted to provide a statistical justification for the minimum $\beta$-aberration criterion, they only obtained a partial relationship between $\beta$-wordlength patterns and contaminations for $k\leq r$.
This result is not sufficient to declare that the minimum $\beta$-aberration criterion minimizes contamination of $k$th-order effects on the estimation of linear effects for all the $k\geq 2$.

In this paper, we study the connection between the minimum contamination criterion and the minimum $\beta$-aberration criterion for the commonly used three-level orthogonal array. 
Through the indicator function, we find that the contamination 
cannot be expressed by a linear combination of the $\beta$-wordlength pattern. 
This result reveals that the minimum $\beta$-aberration criterion is not necessarily equivalent to the minimum contamination criterion for ranking designs. 
%
%
%
We prove that they are theoretically equivalent only when $m=r+1$. 
%
%

At the end of this paper, we would like to emphasize that the minimum $\beta$-aberration design does not minimize contamination of nonnegligible $k$th-order effects on the estimation of linear effects, but on the estimation of the general mean. 
If screening quantitative factors is the main purpose of an experiment, 
the minimum contamination criterion should be more appropriate than the minimum $\beta$-aberration criterion for selecting optimal designs.

\appendix

\section{Proofs}\label{app}

\subsection{Proof of Lemma~\ref{le:l.pq}}
For three-level designs, $c_0^j(x)=1$, $c_1^j(x)={\sqrt \frac{3}{2}}(x-1)$, and $c_2^j(x)=\sqrt 2(\frac{3}{2}(x-1)^2-1)$. By simple calculation, we obtain $c_0^j(x)c_{t_j}^j(x)=c_{t_j}^j(x)$, $c_1^j(x)c_1^j(x)=\frac{1}{\sqrt 2}c_2^j(x)+c_0^j(x)$, $c_1^j(x)c_2^j(x)=\frac{1}{\sqrt 2}c_1^j(x)$, and $c_2^j(x)c_2^j(x)=-\frac{1}{\sqrt 2}c_2^j(x)+c_0^j(x)$. 
Let ${\cal T}_{i,j}^{l(g)}=\{{\bf t}\in{\cal T}\left| \right. t_l=g, ||{\bf t}||_0=i+j, ||{\bf t}||_1=i+2j\}$
and ${\bf 0}|_{u_l=g}$ denote ${\bf u}\in{\cal T}$ with $u_l=g$ and $u_j=0$ for $j=1,\cdots,m,\mbox{ }j\neq l$. 
For given $p$ and $q$, where $p+2q\geq 2$, 
%
\begin{equation}\label{eq:pv.all}
\begin{array}{r@{=}l}
\lambda_{p,q}&\frac{1}{n^2}tr({\bf Z}_{p,q}^T{\bf Z}_1{\bf Z}_1^T{\bf Z}_{p,q})\\
&\frac{1}{n^2}\sum_{l=1}^m\sum_{{\bf t}\in{\cal T}_{p,q}}\left\{{\bf X}_{{\bf 0}|_{u_l=1}}^T{\bf X}_{\bf t}\right\}^2\\
&\frac{1}{n^2}\sum_{l=1}^m\sum_{{\bf t}\in{\cal T}_{p,q}}\left\{\sum_{i=1}^nC_{{\bf 0}|_{u_l=1}}({\bf x}_i) C_{{\bf t}}({\bf x}_i)\right\}^2\\
&\frac{1}{n^2}\sum_{l=1}^m\sum_{{\bf t}\in{\cal T}_{p,q}}\left\{\sum_{i=1}^n\prod_{\substack{ j=1\\j \neq l}}^m[c_{0}^j(x_{ij})c_{t_j}^j(x_{ij})][c_{1}^l(x_{il})c_{t_l}^l(x_{il})]\right\}^2\\
&\frac{1}{n^2}\sum_{l=1}^m\sum_{{\bf t}\in{\cal T}_{p,q}}\left\{\sum_{i=1}^n\prod_{\substack{ j=1\\j \neq l}}^mc_{t_j}^j(x_{ij})[c_{1}^l(x_{il})c_{t_l}^l(x_{il})]\right\}^2.\\
\end{array}
\end{equation}
Because $\sum_{{\bf t}\in{\cal T}_{p,q}}=\sum_{{\bf t}\in{\cal T}_{p,q}^{l(0)}}+\sum_{{\bf t}\in{\cal T}_{p,q}^{l(1)}}+\sum_{{\bf t}\in{\cal T}_{p,q}^{l(2)}}$, we divide equation (\ref{eq:pv.all}) into three parts to discuss.

\noindent Part 1: ${\bf t}\in{\cal T}_{p,q}^{l(0)}$. Because $c_0^l(x_{il})c_{t_l}^l(x_{il})=c_{t_l}^l(x_{il})$, we obtain
\begin{equation}\label{eq:pv.0}
\begin{array}{l}
\frac{1}{n^2}\sum_{l=1}^m\sum_{{\bf t}\in{\cal T}_{p,q}^{l(0)}}\left\{\sum_{i=1}^n\prod_{\substack{ j=1\\j \neq l}}^mc_{t_j}^j(x_{ij})[c_{1}^l(x_{il})c_{0}^l(x_{il})]\right\}^2\\
=\frac{1}{n^2}\sum_{l=1}^m\sum_{{\bf t}'\in{\cal T}_{p+1,q}^{l(1)}}\left\{\sum_{i=1}^n\prod_{\substack{ j=1\\j \neq l}}^mc_{t_j'}^j(x_{ij})c_{1}^l(x_{il})\right\}^2\\
=\sum_{l=1}^m\sum_{{\bf t}'\in{\cal T}_{p+1,q}^{l(1)}}\left (b_{{\bf t}'}/b_{\bf 0}\right )^2.
\end{array}
\end{equation}
For a given $l$, there are ${m-1\choose p}{m-1-p\choose q}$ ${\bf t}'$s in ${\cal T}_{p+1,q}^{l(1)}$. Hence, for $l$ from 1 to $m$, there are total $m {m-1\choose p}{m-1-p\choose q}$ ${\bf t}'$s, where ${\bf t}'$s $\in{\cal T}_{p+1,q}$. Among all the ${\bf t}'$s, there are only ${m\choose p+1}{m-(p+1)\choose q}$ different ${\bf t}'$s $\in{\cal T}_{p+1,q}$ and they appear equally often. Therefore, equation (\ref{eq:pv.0}) equals
\begin{equation}\label{eq:pv.0-1}
\frac{m{m-1\choose p}{m-1-p\choose q}}{{m\choose p+1}{m-(p+1)\choose q}}\sum_{{\bf t}'\in{\cal T}_{p+1,q}}\left (b_{{\bf t}'}/b_{\bf 0}\right )^2
=(p+1)\beta_{p+1,q}.
\end{equation}

\noindent Part 2: ${\bf t}\in{\cal T}_{p,q}^{l(2)}$. Because $c_1^l(x_{il})c_2^l(x_{il})=\frac{1}{\sqrt 2}c_1^l(x_{il})$, we obtain
\begin{equation}\label{eq:pv.2}
\begin{array}{l}
\frac{1}{n^2}\sum_{l=1}^m\sum_{{\bf t}\in{\cal T}_{p,q}^{l(2)}}\left\{\sum_{i=1}^n\prod_{\substack{ j=1\\j \neq l}}^mc_{t_j}^j(x_{ij})[c_{1}^l(x_{il})c_{2}^l(x_{il})]\right\}^2\\
=\frac{1}{n^2}\sum_{l=1}^m\sum_{{\bf t}'\in{\cal T}_{p+1,q-1}^{l(1)}}\left\{\sum_{i=1}^n\prod_{\substack{ j=1\\j \neq l}}^mc_{t_j'}^j(x_{ij})[\frac{1}{\sqrt 2}c_{1}^l(x_{il})]\right\}^2\\
=\frac{1}{2}\sum_{l=1}^m\sum_{{\bf t}'\in{\cal T}_{p+1,q-1}^{l(1)}}\left (b_{{\bf t}'}/b_{\bf 0}\right )^2.
\end{array}
\end{equation}
For a given $l$, there are ${m-1\choose p}{m-1-p\choose q-1}$ ${\bf t}'$s in ${\cal T}_{p+1,q-1}^{l(1)}$. Hence, for $l$ from 1 to $m$, there are total $m {m-1\choose p}{m-1-p\choose q-1}$ ${\bf t}'$s, where ${\bf t}'$s $\in{\cal T}_{p+1,q-1}$. Among all the ${\bf t}'$s, there are only ${m\choose p+1}{m-(p+1)\choose q-1}$ different ${\bf t}'$s $\in{\cal T}_{p+1,q-1}$ and they appear equally often. Therefore, equation (\ref{eq:pv.2}) equals
\begin{equation}\label{eq:pv.2-1}
\frac{m{m-1\choose p}{m-1-p\choose q-1}}{2{m\choose p+1}{m-(p+1)\choose q-1}}\sum_{{\bf t}'\in{\cal T}_{p+1,q-1}}\left (b_{{\bf t}'}/b_{\bf 0}\right )^2
=\frac{p+1}{2}\beta_{p+1,q-1}.
\end{equation}

\noindent Part 3: ${\bf t}\in{\cal T}_{p,q}^{l(1)}$. Because $c_1^l(x_{il})c_1^l(x_{il})=\frac{1}{\sqrt 2}c_2^l(x_{il})+c_0^l(x_{il})$, we obtain
\begin{equation}\label{eq:pv.1}
\begin{array}{l}
\frac{1}{n^2}\sum_{l=1}^m\sum_{{\bf t}\in{\cal T}_{p,q}^{l(1)}}\left\{\sum_{i=1}^n\prod_{\substack{ j=1\\j \neq l}}^mc_{t_j}^j(x_{ij})[c_{1}^l(x_{il})c_{1}^l(x_{il})]\right\}^2\\
=\frac{1}{n^2}\sum_{l=1}^m\sum_{{\bf t}'\in{\cal T}_{p-1,q+1}^{l(2)}}\left\{\sum_{i=1}^n\prod_{\substack{ j=1\\j \neq l}}^mc_{t_j'}^j(x_{ij})[\frac{1}{{\sqrt 2}}c_{2}^l(x_{il})]\right\}^2\\
\mbox{ }\mbox{ }\mbox{ }+\frac{1}{n^2}\sum_{l=1}^m\sum_{{\bf t}''\in{\cal T}_{p-1,q}^{l(0)}}\left\{\sum_{i=1}^n\prod_{\substack{ j=1\\j \neq l}}^mc_{t_j''}^j(x_{ij})c_{0}^l(x_{il})\right\}^2\\
\mbox{ }\mbox{ }\mbox{ }+\frac{2}{n^2}\sum_{l=1}^m\sum_{{\bf t}^*\in{\cal T}_{p-1,q}^{l(-)}}\left\{[\frac{1}{{\sqrt 2}}\sum_{i=1}^n\prod_{\substack{ j=1\\j \neq l}}^mc_{t_j^*}^j(x_{ij})c_{2}^l(x_{il})]\right.\\
\mbox{ }\mbox{ }\mbox{ }\mbox{ }\mbox{ }\mbox{ }\mbox{ }\mbox{ }\mbox{ }\mbox{ }\mbox{ }\mbox{ }\mbox{ }\mbox{ }\mbox{ }\mbox{ }\mbox{ }\mbox{ }\mbox{ }\mbox{ }\mbox{ }\mbox{ }\mbox{ }\mbox{ }\mbox{ }\mbox{ }\mbox{ }\mbox{ }\mbox{ }\mbox{ }\mbox{ }\mbox{ }\mbox{ }\mbox{ }\mbox{ }\mbox{ }\times\left.[\sum_{i=1}^n\prod_{\substack{ j=1\\j \neq l}}^mc_{t_j^*}^j(x_{ij})c_{0}^l(x_{il})]\right\}\\
=\frac{1}{2}\sum_{l=1}^m\sum_{{\bf t}'\in{\cal T}_{p-1,q+1}^{l(2)}}\left (b_{{\bf t}'}/b_{\bf 0}\right )^2
+\sum_{l=1}^m\sum_{{\bf t}''\in{\cal T}_{p-1,q}^{l(0)}}\left (b_{{\bf t}''}/b_{\bf 0}\right )^2\\
\mbox{ }\mbox{ }\mbox{ }+{\sqrt 2}\sum_{l=1}^m\sum_{{\bf t}^*\in{\cal T}_{p-1,q}^{l(-)}}\prod_{g=0,2}(b_{{\bf t}^*|_{t_l^*=g}}/b_{\bf 0}).
\end{array}
\end{equation}
For a given $l$, there are ${m-1\choose p-1}{m-1-(p-1)\choose q}$ ${\bf t}'$s in ${\cal T}_{p-1,q+1}^{l(2)}$. Hence, for $l$ from 1 to $m$, there are total $m{m-1\choose p-1}{m-1-(p-1)\choose q}$ ${\bf t}'$s, where ${\bf t}'$s $\in{\cal T}_{p-1,q+1}$. Among all the ${\bf t}'$s, there are only ${m\choose p-1}{m-(p-1)\choose q+1}$ different ${\bf t}'$s $\in{\cal T}_{p-1,q+1}$ and they appear equally often. 
Similarly, for a given $l$, there are ${m-1\choose p-1}{m-1-(p-1)\choose q}$ ${\bf t}''$s in ${\cal T}_{p-1,q}^{l(0)}$. Hence, for $l$ from 1 to $m$, there are total $m{m-1\choose p-1}{m-1-(p-1)\choose q}$ ${\bf t}''$s, where ${\bf t}''$s $\in{\cal T}_{p-1,q}$. Among all the ${\bf t}''$s, there are only ${m\choose p-1}{m-(p-1)\choose q}$ different ${\bf t}''$s $\in{\cal T}_{p-1,q}$ and they appear equally often. Therefore, equation (\ref{eq:pv.1}) equals
\begin{equation}\label{eq:pv.1-1}
\begin{array}{l}
\frac{m{m-1\choose p-1}{m-1-(p-1)\choose q}}{2{m\choose p-1}{m-(p-1)\choose q+1}}\sum_{{\bf t}'\in{\cal T}_{p-1,q+1}}\left (b_{{\bf t}'}/b_{\bf 0}\right )^2
+\frac{m{m-1\choose p-1}{m-1-(p-1)\choose q}}{{m\choose p-1}{m-(p-1)\choose q}}\sum_{{\bf t}''\in{\cal T}_{p-1,q}}\left (b_{{\bf t}''}/b_{\bf 0}\right )^2\\
\mbox{ }\mbox{ }\mbox{ }+{\sqrt 2}\sum_{l=1}^m\sum_{{\bf t}^*\in{\cal T}_{p-1,q}^{l(-)}}\prod_{g=0,2}(b_{{\bf t}^*|_{t_l^*=g}}/b_{\bf 0})\\
=\frac{q+1}{2}\beta_{p-1,q+1}
+(m-p-q+1)\beta_{p-1,q}
+{\sqrt 2}\xi_{p-1,q}.
\end{array}
\end{equation}
The result holds by summarizing equations (\ref{eq:pv.0-1}), (\ref{eq:pv.2-1}), and (\ref{eq:pv.1-1}).



\subsection{Proof of Proposition~\ref{pr:l.k}} 
Recall that $\beta_{i,j}=0$ for $i,j<0$, $i,j>m$, or $i+j>m$, and $\xi_{i,j}=0$ if $i,j<0$, $i,j>m-1$, or $i+j>m-1$. When $k$ is even, according to equations (\ref{eq:b.k=b.pq}) and (\ref{eq:l.k=l.pq}) and Lemma~\ref{le:l.pq}, we obtain 
\begin{equation}\label{eq:pv.even}
\begin{array}{rcl}
\lambda_k&=&\sum_{j'=0}^{k/2}\lambda_{k-2j',j'}\\
&=&\left\{\sum_{j'=0}^{k/2}[(k-2j')+1]\beta_{(k-2j')+1,j'}+\sum_{j'=0}^{k/2-1}\frac{j'+1}{2}\beta_{(k-2j')-1,j'+1}\right\}\\
&&+\left\{\sum_{j'=1}^{k/2}\frac{(k-2j')+1}{2}\beta_{(k-2j')+1,j'-1}+\sum_{j'=0}^{k/2-1}(m-(k-2j')-j'+1)\beta_{(k-2j')-1,j'}\right\}\\
&&+\sqrt2\sum_{j'=0}^{k/2-1}\xi_{(k-2j')-1,j'}\\
&=&\left\{\sum_{j=0}^{k/2}(k-2j+1)\beta_{k-2j+1,j}+\sum_{j=1}^{k/2}\frac{j}{2}
\beta_{k-2j+1,j}\right\}\\
&&+\left\{\sum_{j=0}^{k/2-1}
\frac{k-2j-1}{2}\beta_{k-2j-1,j}
+\sum_{j=0}^{k/2-1}(m-k+j+1)\beta_{k-2j-1,j}\right\}\\
&&+\sqrt2\sum_{j=0}^{k/2-1}\xi_{k-2j-1,j}\\
&=&\sum_{j=0}^{k/2}(k-\frac{3}{2}j+1)\beta_{k-2j+1,j}+(m-\frac{k-1}{2})\sum_{j=0}^{k/2-1}\beta_{k-2j-1,j}\\
&&+\sqrt2\sum_{j=0}^{k/2-1}\xi_{k-2j-1,j}\\
&=&(k-\frac{3}{2}\times\frac{k}{2}+1)\sum_{j=0}^{k/2}\beta_{k-2j+1,j}+\sum_{j=0}^{k/2-1}(\frac{3}{2}\times\frac{k}{2}-\frac{3}{2}j)
\beta_{k-2j+1,j}\\
&&+(m-\frac{k-1}{2})\sum_{j=0}^{k/2-1}\beta_{k-2j-1,j}+\sqrt2\sum_{j=0}^{k/2-1}\xi_{k-2j-1,j}\\
&=&(k-\frac{3}{2}\times\frac{k}{2}+1)\sum_{j=0}^{\lfloor(k+1)/2\rfloor}\beta_{(k+1)-2j,j}+(m-\frac{k-1}{2})\sum_{j=0}^{\lfloor (k-1)/2 \rfloor}\beta_{(k-1)-2j,j}\\
&&+\frac{3}{2}\sum_{j=0}^{k/2-1}(\frac{k}{2}-j)\beta_{k-2j+1,j}+\sqrt2\sum_{j=0}^{k/2-1}\xi_{k-2j-1,j}\\
&=&(1+k-\frac{3}{2}\times\lceil\frac{k}{2}\rceil)\beta_{k+1}+(m-\frac{k-1}{2})\beta_{k-1}\\
&&+\frac{3}{2}\sum_{j=0}^{\lceil k/2\rceil-1}(\lceil\frac{k}{2}\rceil-j)\beta_{k-2j+1,j}+{\sqrt 2}\sum_{j=0}^{\lceil k/2\rceil-1} \xi_{k-2j-1,j}.
\end{array}
\end{equation}
When $k$ is odd, we obtain
\begin{equation}\label{eq:pv.odd}
\begin{array}{rcl}
\lambda_k&=&\sum_{j'=0}^{(k-1)/2}\lambda_{k-2j',j'}\\
&=&\left\{\sum_{j'=0}^{(k-1)/2}[(k-2j')+1]\beta_{(k-2j')+1,j'}+\sum_{j'=0}^{(k-1)/2}\frac{j'+1}{2}\beta_{(k-2j')-1,j'+1}\right\}\\
&&+\left\{\sum_{j'=1}^{(k-1)/2}\frac{(k-2j')+1}{2}\beta_{(k-2j')+1,j'-1}+\sum_{j'=0}^{(k-1)/2}(m-(k-2j')-j'+1)\beta_{(k-2j')-1,j'}\right\}\\
&&+\sqrt2\sum_{j'=0}^{(k-1)/2}\xi_{(k-2j')-1,j'}\\
&=&\left\{\sum_{j=0}^{(k-1)/2}(k-2j+1)\beta_{k-2j+1,j}+\sum_{j=1}^{(k-1)/2+1}\frac{j}{2}
\beta_{k-2j+1,j}\right\}\\
&&+\left\{\sum_{j=0}^{(k-1)/2-1}
\frac{k-2j-1}{2}\beta_{k-2j-1,j}
+\sum_{j=0}^{(k-1)/2}(m-k+j+1)\beta_{k-2j-1,j}\right\}\\
&&+\sqrt2\sum_{j=0}^{(k-1)/2}\xi_{k-2j-1,j}\\
&=&\sum_{j=0}^{(k+1)/2}(k-\frac{3}{2}j+1)\beta_{k-2j+1,j}+(m-\frac{k-1}{2})\sum_{j=0}^{(k-1)/2}\beta_{k-2j-1,j}\\
&&+\sqrt2\sum_{j=0}^{(k-1)/2}\xi_{k-2j-1,j}\\
&=&(k-\frac{3}{2}\times\frac{k+1}{2}+1)\sum_{j=0}^{(k+1)/2}\beta_{k-2j+1,j}+\sum_{j=0}^{(k-1)/2}(\frac{3}{2}\times\frac{k+1}{2}-\frac{3}{2}j)\beta_{k-2j+1,j}\\
&&+(m-\frac{k-1}{2})\sum_{j=0}^{(k-1)/2}\beta_{k-2j-1,j}+\sqrt2\sum_{j=0}^{(k-1)/2}\xi_{k-2j-1,j}\\
&=&(k-\frac{3}{2}\times\frac{k+1}{2}+1)\sum_{j=0}^{\lfloor (k+1)/2\rfloor}\beta_{(k+1)-2j,j}+(m-\frac{k-1}{2})\sum_{j=0}^{\lfloor (k-1)/2\rfloor}\beta_{(k-1)-2j,j}\\
&&+\frac{3}{2}\sum_{j=0}^{(k-1)/2}(\frac{k+1}{2}-j)\beta_{k-2j+1,j}+\sqrt2\sum_{j=0}^{(k-1)/2}\xi_{k-2j-1,j}\\
&=&(1+k-\frac{3}{2}\times\lceil\frac{k}{2}\rceil)\beta_{k+1}+(m-\frac{k-1}{2})\beta_{k-1}\\
&&+\frac{3}{2}\sum_{j=0}^{\lceil k/2\rceil-1}(\lceil\frac{k}{2}\rceil-j)\beta_{k-2j+1,j}+{\sqrt 2}\sum_{j=0}^{\lceil k/2\rceil-1} \xi_{k-2j-1,j}.
\end{array}
\end{equation}

\subsection{Proof of Proposition~\ref{pr:r+1}}
We need the following Lemma to prove Proposition~\ref{pr:r+1}.


\begin{lemma}\label{le:=0} 
If $0<i+j\leq r$, then (a) $\beta_{i,j}=0$ and (b) $\xi_{i,j}=0$.
\end{lemma}

\textbf{Proof of Lemma~\ref{le:=0}.}
For orthogonal arrays with strength $r$, $b_{\bf t}=0$ if $0<||{\bf t}||_0\leq r$. 
Because $\beta_{i,j}$ is the summation of $(b_{\bf t}/b_{\bf 0})^2$ for all $\bf t$s such that $||{\bf t}||_0=i+j$, it is obvious that  
result (a) holds if $0< i+j\leq r$.
For ${\bf t}\in{\cal T}_{i,j}^{l(-)}$, $||{\bf t}|_{t_l=0}||_0=i+j$ and hence $b_{{\bf t}|_{t_l=0}}=0$ if $0<i+j\leq r$. 
We obtain $\xi_{i,j}=\sum_{l=1}^m\sum_{{\bf t}\in{\cal T}_{i,j}^{l(-)}}(0/b_{\bf 0})\times(b_{{\bf t}|_{t_l=2}}/b_{\bf 0})=0$ when $0<i+j\leq r$.\\

\textbf{Proof of Proposition~\ref{pr:r+1}.}
Let us rewrite equation (\ref{eq:l.k}) as 
\begin{equation}\label{eq:re}
\lambda_{k}=\left\{(1+k-\frac{3}{2}\times\lceil \frac{k}{2}\rceil)\beta_{k+1}+B_1+B_2\right\}+(m-\frac{k-1}{2})\beta_{k-1},
\end{equation}
where $B_1=\frac{3}{2}\sum_{j=0}^{\lceil k/2\rceil-1}(\lceil\frac{k}{2}\rceil-j)\beta_{k-2j+1,j}$ and $B_2=
\sqrt 2 \sum_{j=0}^{\lceil k/2\rceil-1}\xi_{k-2j-1,j}$. We first prove $B_2=0$ when $m=r+1$. For convenience, let $i=k-2j-1$. 
%
For $k=2,\cdots,m$, $B_2=\sqrt 2 \sum_{j=0}^{\lceil k/2\rceil-1}\xi_{i,j}$. However, $\xi_{i,j}=0$ since $i+j=k-1-j\leq k-1\leq m-1=r$ (Lemma~\ref{le:=0} (b)).
%
For $k=m+1,\cdots,m'$, $\xi_{i,j}=0$ for $j=0,\cdots,k-m-1$ since $i+j=k-1-j\geq k-1-(k-m-1)=m$ (according to the definition). Hence, $B_2=\sqrt 2 \sum_{j=k-m}^{\lceil k/2\rceil-1}\xi_{i,j}$. Since $i+j\leq k-1-(k-m)=m-1=r$ (Lemma~\ref{le:=0} (b)), $\xi_{i,j}=0$.
Therefore, $B_2=0$. 
Now let us prove $(1+k-\frac{3}{2}\times\lceil \frac{k}{2}\rceil)\beta_{k+1}+B_1=\rho\beta_{k+1}$ when $m=r+1$, where $\rho$ is defined in equation (\ref{eq:rho}).
%
%
For $k=2,\cdots,r$, $\beta_{k-2j+1,j}=0$ for $j=1,\cdots,\lceil k/2\rceil$ since $(k-2j+1)+j=k+1-j\leq k+1-(1)=k\leq r$ (Lemma~\ref{le:=0}~(a)). We obtain
\begin{equation}\label{eq:ap1}
\begin{array}{l}
(1+k-\frac{3}{2}\lceil\frac{k}{2}\rceil)\beta_{k+1}+B_1\\
=\sum_{j=0}^{\lceil k/2\rceil}(k-\frac{3}{2}j+1)\beta_{k-2j+1,j}\\
=(k+1)\beta_{k+1,0}\\
=(k+1)\sum_{j'=0}^{\lfloor (k+1)/2\rfloor}\beta_{(k+1)-2j',j'}\\
=(k+1)\beta_{k+1}.
\end{array}
\end{equation}
For $k=r+1,\cdots,m'-1$, $\beta_{k-2j+1,j}=0$ for $j=0,\cdots,k-r-1$ since $(k-2j+1)+j=k+1-j\geq k+1-(k-r-1)=r+2=m+1$ (according to the definition) and $\beta_{k-2j+1,j}=0$ for $j=k-r+1,\cdots,\lceil k/2\rceil$ since $(k-2j+1)+j=k+1-j\leq k+1-(k+1-r)=r$ (Lemma~\ref{le:=0}~(a)). 
It follows that
\begin{equation}\label{eq:ap2}
\begin{array}{l}
(1+k-\frac{3}{2}\lceil\frac{k}{2}\rceil)\beta_{k+1}+B_1\\
=\sum_{j=k-r}^{\lceil k/2\rceil}(k-\frac{3}{2}j+1)\beta_{k-2j+1,j}\\
=(k-\frac{3}{2}(k-r)+1)\beta_{2r-k+1,k-r}\\
=\frac{1}{2}(3r+2-k)\beta_{2r-k+1,k-r}\\
=\frac{1}{2}(3r+2-k)\sum_{j'=0}^{\lfloor (k+1)/2\rfloor}\beta_{(k+1)-2j',j'}\\
=\frac{1}{2}(3r+2-k)\beta_{k+1}.
\end{array}
\end{equation}
For $k=m' (=2m)$, $\beta_{k-2j+1,j}=0$ for $j=0,\cdots,\lceil k/2\rceil$ since $(k-2j+1)+j=k+1-j\geq k+1-\lceil k/2\rceil=2m+1-m=m+1$ (according to the definition). Therefore, 
\begin{equation}\label{eq:ap3}
\begin{array}{l}
(1+k-\frac{3}{2}\lceil\frac{k}{2}\rceil)\beta_{k+1}+B_1\\
=\sum_{j=0}^{\lceil k/2\rceil}(k-\frac{3}{2}j+1)\beta_{k-2j+1,j}\\
=0.
\end{array}
\end{equation} 
The result holds by replacing equations (\ref{eq:ap1}), (\ref{eq:ap2}), (\ref{eq:ap3}) and $m=r+1$ into equation (\ref{eq:re}).

\subsection{Proof of Theorem~\ref{th:eq}} The result follows according to Proposition~\ref{pr:r+1}.

\subsection{Proof of Theorem~\ref{th:even}}
According to Theorem 2 in \citet{r12}, 
for a mirror-symmetric design, we have $\beta_{k}=0$ if $k$ is odd, which also implies $\beta_{i,j}=0$ and $b_{{\bf t}|_{t_l=0}}=0$ for ${\bf t}\in{\cal T}_{i,j}^{l(-)}$ if $i+2j$ is odd. 
Now let us consider $\lambda_k$ when $k$ is even. 
Recall that 
\begin{equation}\label{eq:ap.re}
\begin{array}{rl}
\lambda_{k}=&(1+k-\frac{3}{2}\times\lceil \frac{k}{2}\rceil)\beta_{k+1}\\
&+(m-\frac{k-1}{2})\beta_{k-1}\\
&+\frac{3}{2}\sum_{j=0}^{\lceil k/2\rceil-1}(\lceil\frac{k}{2}\rceil-j)\beta_{k-2j+1,j}\\
&+\sqrt 2 \sum_{j=0}^{\lceil k/2\rceil-1}\xi_{k-2j-1,j}.
\end{array}
\end{equation}
We will prove that the four terms in equation~(\ref{eq:ap.re}) are all zero when $k$ is even. 
First, $\beta_{k+1}=0$ and $\beta_{k-1}=0$ since $k+1$ and $k-1$ is odd. Hence, $(1+k-\frac{3}{2}\times\lceil \frac{k}{2}\rceil)\beta_{k+1}=0$ and $(m-\frac{k-1}{2})\beta_{k-1}=0$. 
Second, $\beta_{k-2j+1,j}=0$ since $(k-2j+1)+2j=k+1$ (odd). Hence, $\frac{3}{2}\sum_{j=0}^{\lceil k/2\rceil-1}(\lceil\frac{k}{2}\rceil-j)\beta_{k-2j+1,j}=0$. 
Finally, $b_{{\bf t}|_{t_l=0}}=0$ for ${\bf t}\in{\cal T}_{k-2j-1,j}^{l(-)}$ since $(k-2j-1)+2j=k-1$ (odd). Hence, $\sqrt 2 \sum_{j=0}^{\lceil k/2\rceil-1}\xi_{k-2j-1,j}=\sqrt 2 \sum_{j=0}^{\lceil k/2\rceil-1}\sum_{l=1}^m\sum_{{\bf t}\in{\cal T}_{k-2j-1,j}^{l(-)}}(0/b_{\bf 0})\times(b_{{\bf t}|_{t_l=2}}/b_{\bf 0})=0$. 
The result follows by summarizing above discussions.




\bibliographystyle{imsart-nameyear}
\bibliography{2014-beta}

\begin{thebibliography}{6}

\bibitem[\protect\citeauthoryear{Cheng and Ye}{2004}]{r1}
\begin{barticle}[author]
\bauthor{\bsnm{Cheng},~\bfnm{Shao-Wei}\binits{S.-W.}} \AND
  \bauthor{\bsnm{Ye},~\bfnm{Kenny~Q.}\binits{K.~Q.}}
(\byear{2004}).
\btitle{Geometric isomorphism and minimum aberration for factorial designs with
  quantitative factors}.
\bjournal{Ann. Statist.}
\bvolume{32}
\bpages{2168--2185}.
\bmrnumber{2102507}
\end{barticle}
\endbibitem

\bibitem[\protect\citeauthoryear{Fontana, Pistone and Rogantin}{2000}]{r4}
\begin{barticle}[author]
\bauthor{\bsnm{Fontana},~\bfnm{Roberto}\binits{R.}},
  \bauthor{\bsnm{Pistone},~\bfnm{Giovanni}\binits{G.}} \AND
  \bauthor{\bsnm{Rogantin},~\bfnm{Maria~Piera}\binits{M.~P.}}
(\byear{2000}).
\btitle{Classification of two-level factorial fractions}.
\bjournal{J. Statist. Plann. Inference}
\bvolume{87}
\bpages{149--172}.
\bmrnumber{1772046}
\end{barticle}
\endbibitem

\bibitem[\protect\citeauthoryear{Tang and Deng}{1999}]{r11}
\begin{barticle}[author]
\bauthor{\bsnm{Tang},~\bfnm{Boxin}\binits{B.}} \AND
  \bauthor{\bsnm{Deng},~\bfnm{Lih-Yuan}\binits{L.-Y.}}
(\byear{1999}).
\btitle{Minimum {$G\sb 2$}-aberration for nonregular fractional factorial
  designs}.
\bjournal{Ann. Statist.}
\bvolume{27}
\bpages{1914--1926}.
\bmrnumber{1765622}
\end{barticle}
\endbibitem

\bibitem[\protect\citeauthoryear{Tang and Xu}{2014}]{r12}
\begin{barticle}[author]
\bauthor{\bsnm{Tang},~\bfnm{Yu}\binits{Y.}} \AND
  \bauthor{\bsnm{Xu},~\bfnm{Hongquan}\binits{H.}}
(\byear{2014}).
\btitle{Permuting regular fractional factorial designs for screening
  quantitative factors}.
\bjournal{Biometrika}
\bvolume{101}
\bpages{333--350}.
\end{barticle}
\endbibitem

\bibitem[\protect\citeauthoryear{Xu and Wu}{2001}]{r14}
\begin{barticle}[author]
\bauthor{\bsnm{Xu},~\bfnm{Hongquan}\binits{H.}} \AND
  \bauthor{\bsnm{Wu},~\bfnm{C.~F.~J.}\binits{C.~F.~J.}}
(\byear{2001}).
\btitle{Generalized minimum aberration for asymmetrical fractional factorial
  designs}.
\bjournal{Ann. Statist.}
\bvolume{29}
\bpages{1066--1077}.
\bmrnumber{1863969}
\end{barticle}
\endbibitem

\bibitem[\protect\citeauthoryear{Ye}{2003}]{r15}
\begin{barticle}[author]
\bauthor{\bsnm{Ye},~\bfnm{Kenny~Q.}\binits{K.~Q.}}
(\byear{2003}).
\btitle{Indicator function and its application in two-level factorial designs}.
\bjournal{Ann. Statist.}
\bvolume{31}
\bpages{984--994}.
\bmrnumber{1994738}
\end{barticle}
\endbibitem

\end{thebibliography}


\begin{thebibliography}{9}

\bibitem{r1}
\textsc{Cheng, S.-W.} and \textsc{Ye, K. Q.} (2004). 
Geometric isomorphism and minimum aberration for factorial designs with quantitative factors.
\textit{Ann. Statist.}
\textbf{32}, 2168--2185.
\MR{2102507}

\bibitem{r2}
\textsc{Fontana, R.} and \textsc{Pistone, G.} and \textsc{Rogantin, M. P.} (2000). 
Classification of two-level factorial fractions.
\textit{J. Statist. Plann. Inference.}
\textbf{87}, 149--172.
\MR{1772046}

\bibitem{r3}
\textsc{Tang, B.} and \textsc{Deng, L.-Y.} (1999). 
Minimum $G\sb 2$-aberration for nonregular fractional factorial designs.
\textit{Ann. Statist.}
\textbf{27} 1914--1926.
\MR{1765622}

\bibitem{r4}
\textsc{Tang, Y.} and \textsc{Xu, H.}  (2014). 
Permuting regular fractional factorial designs for screening quantitative factors.
\textit{Biometrika}
\textbf{101} 333--350.
\MR{}

\bibitem{r5}
\textsc{Xu, H.} and \textsc{Wu, C.F.J.}  (2001). 
Generalized minimum aberration for asymmetrical fractional factorial designs.
\textit{Ann. Statist.}
\textbf{29} 1066--1077.
\MR{1869240}

\bibitem{r6}
\textsc{Ye, K. Q.} (2003).
Indicator function and its application in two-level factorial designs.
\textit{Ann. Statist.}
\textbf{31} 984--994.
\MR{1994738}


\end{thebibliography}

\end{document}